\documentclass[a4paper]{article}

\usepackage{epsfig}

\usepackage{amssymb}
\title{Coding rotations on intervals}
\author{Jean Berstel\\
\footnotesize{Institut Gaspard Monge (IGM)}\\[-6pt]
\footnotesize{Universit\'e de Marne-la-Vall\'ee}\\[-6pt]
\footnotesize{5, boulevard Descartes, 77454 Marne-la-Vall\'ee C\'{e}dex 2}\\
\and Laurent Vuillon\\
\footnotesize{Laboratoire d'informatique algorithmique: fondements  et
applications (LIAFA)}\\[-6pt]
\footnotesize{Universit\'e Denis-Diderot (Paris VII)}\\[-6pt]
\footnotesize{2, place Jussieu, 75251 Paris C\'{e}dex 05}
\date{}
}
\def\petitcarre{\vrule height4pt width 4pt depth0pt}
\def\QED{\relax\ifmmode\eqno{\hbox{\petitcarre}}\else
{\unskip\nobreak\hfil\penalty50
 \hskip2em\hbox{}\nobreak\hfil\petitcarre 
 \parfillskip=0pt \finalhyphendemerits=0\par\medskip\par}\fi}

\newenvironment{proof}{\medskip\par\noindent\textit{Proof}}{\QED}
\newtheorem{theorem}{Theorem}[section]
\newtheorem{proposition}[theorem]{Proposition}

\newtheorem{lemma}[theorem]{Lemma}

\def\A{{\cal A}}

\newcommand{\nat}{\mathbb{N}}

\newcommand{\Z}{\mathbb{Z}}
\newcommand{\R}{\mathbb{R}}
\newcommand{\T}{\mathbb{T}}

\begin{document}
\maketitle

\begin{abstract}\noindent
 We show that the coding of a rotation by 
$\alpha$ on $m$ intervals with rationally
independent
lengths can be recoded over $m$ Sturmian words of angle $\alpha$.
\end{abstract}


\baselineskip=1.3\baselineskip

\parskip=\medskipamount

\section{Introduction}

The coding of rotations is a tool for the construction of infinite words over
a finite alphabet. Consider
a rotation $R_\alpha$, given by an angle $\alpha$,
and defined for a point $x$ by
$R_\alpha(x)=\{x+\alpha\}$ where $\{y\}$
denotes the fractional part of $y$.
Consider next a partition of the unit circle in $m$
 half open intervals $\{I_1,I_2, \cdots I_m\}$. For any 
starting point $x$ with $0 \leq x < 1$, one
gets
 an infinite word $u$ by $ I(x) I(R_\alpha(x))
I(R^2_\alpha(x))
 \cdots I(R^n_\alpha(x)) \cdots$, where $I(y)=i$ if  $y \in I_i$.

In the special case where $\alpha$ is irrational and
and the partition  is $I_1=[0, \alpha[$ and $I_2=[\alpha,1[$,
this construction produces exactly the well-known Sturmian words.
These words appear in various domains as computer sciences \cite{BS},
Physics, Mathematical optimization
and play a crucial role in this article. It is remarkable that
Sturmian words have a combinatorial characterization. Thus, they are exactly
aperiodic words with (subword)
complexity $p(n)=n+1$ where the complexity function $p:\nat\to\nat$
counts the number of distinct factors of length $n$ in the
infinite word $u$ \cite{BS}. The same general construction allows
also to compute Rote words with complexity $p(n)=2n$ by using an
irrational rotation
and the partition $I_1=[0, \frac{1}{2}[$ and $I_2=[\frac{1}{2},1[$
(see \cite{Rot}). More generally, one can obtain
infinite words with complexity $p(n)=an+b,$ where $a$ and $b$ are real,
by coding of rotation
 \cite{AB,BerH}.

In addition,  codings of rotation with an irrational
value of $\alpha$
and the partition $I_1=[0, \beta[$ and $I_2=[\beta,1[$ are intimately
related to Sturmian words.
Indeed, the first sequence is the difference term by term of two
Sturmian words \cite{Hub}.
 Didier gives a characterization of the coding of rotation with a
partition of $m$ intervals of length greater than
$\alpha$
by using Sturmian words and  cellular automata \cite{Did}.
Finally, Blanchard and Kurka study the complexity of formal languages
that are generated by coding of rotation \cite{BK}.

\begin{center}
\begin{figure}[htb]
\input{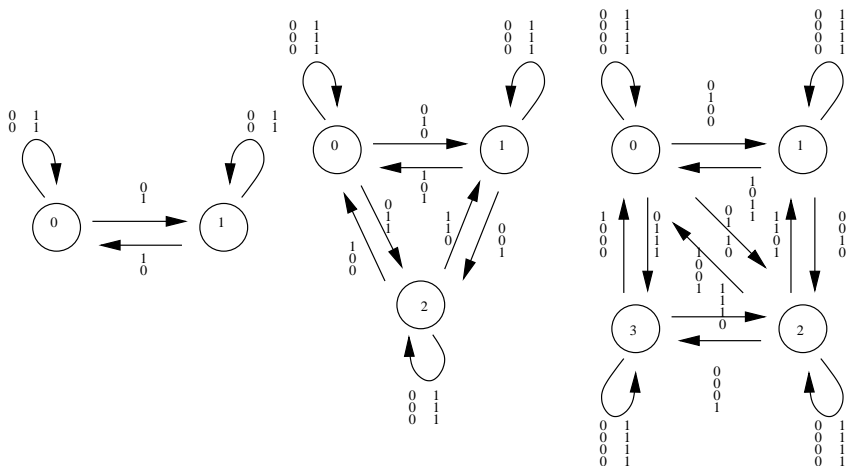}
\caption{Automata for $m=1,2,3$.} 
\label{Tr}
\end{figure}
\end{center}

The goal of this article is to show that the coding of a rotation by 
$\alpha$ on $m$ intervals with rationally
independent
lengths can be recoded over $m$ Sturmian words of angle $\alpha$.
More precisely,   for a given $m$ an universal automaton is
constructed  such that the edge indexed by
the vector
of values of the $i$th letter on each Sturmian word gives the value of
the $i$th letter of the coding
of rotation (see Figure \ref{Tr}). If the partition is given by $[\beta_j,\beta_{j+1}[$ where
$\beta_0 = 0 < \beta_1 < \beta_2 <
 \cdots < \beta_j < \cdots < \beta_{m+1}=1$,
 then the $\ell$th Sturmian word
is given by the partition
$I_1=[\beta_\ell, \beta_\ell + \alpha \bmod1[$ and the
complement of $I_1$ on the unit circle. 

\section{Examples}

\begin{center}
\begin{figure}[htb]
\input{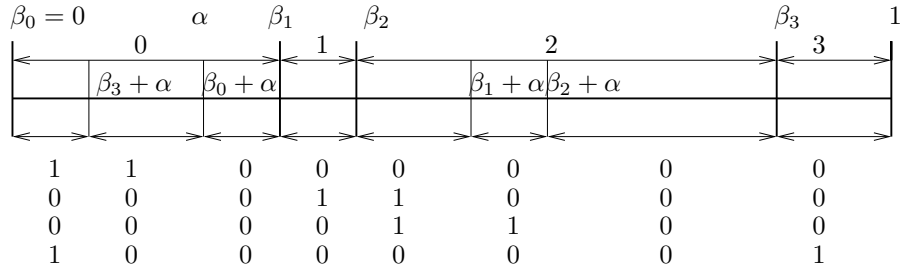}
\caption{Partition of the unit circle.} 
\label{part}
\end{figure}
\end{center}

The figure \ref{part} shows a
partition of the unit circle by 4 intervals of form
$[\beta_j,\beta_{j+1}[$
and the coding by 8 intervals associated with  binary vectors
(we can find the coding of the interval
$[\beta_j,\beta_{j+1}[$) by the automaton for $m=3$ applied to
 the binary vector value.

As an example, using
the universal automaton for $m=2$, the three following Sturmian words can be recoded   on a
word on a three letter alphabet.
$$1001010010100101 \cdots$$
$$0100101001010010 \cdots $$
$$0010100101001010 \cdots  $$
is recoded on the following  word:

$$0120201202012020   \cdots  .$$

\section{Notation}

We will consider subsets of $[0,1[$ that we call intervals. Let $x,y$ be in $[0,1[$. Then we set
$$[x,y[  \cases{
\{z\mid x\le z < y\}& if $x<y$\cr
\emptyset & if $x=y$\cr
\{z\mid x\le z < 1\}\cup \{z\mid 0\le z < y\}& if $x>y$
}$$
In particular, $[x,y[ = [0,y[\,\cup\, [x,1[$ if $x>y$.
This is precisely the notion of an interval on the torus $\T=\R/\Z$. 

Let $\alpha$, $\beta_1,\ldots,\beta_m$ be numbers in the
interval $]0,1[$, with $\beta_1<\cdots<\beta_m$. It will be convenient to set $\beta_0=0$ and $\beta_{m+1}=1$. The $m+1$ intervals
$$B_k=[\beta_k,\beta_{k+1}[,\qquad k=0,\ldots,m$$
are a partition of $[0,1[$. We consider the \textit{rotation} of angle $\alpha$ defined by $R_\alpha(x) = x+\alpha \bmod1$. Define intervals $I_k$ by
(all values are computed modulo $1$)
$$I_k=[\beta_k,\beta_k+\alpha[, \qquad k=0,\ldots,m$$
We will be specially interested in the nonempty intervals
$$X_K = \bigcap_{k\in K} I_k\cap \bigcap_{k\notin K} \overline{I_k}$$
Here, $K$ is a subset of $M=\{0,\ldots, m\}$, and 
$\overline{I_k}=[0,1[\setminus I_k$ is the complement of
$I_k$. Observe that, for any nonempty interval $I=[x,y[$, one gets
$\overline{I}=[y,x[$.

\section{Circular order}

We want to compute intersections of intervals. Although the geometric
approach is easy to understand, it is error prone because points are
usually not in general position. Therefore, we consider a more
combinatoric approach.

Given numbers $x_1,\ldots, x_n\in [0,1[$, the sequence $(x_1,\ldots,
x_n)$ is \textit{circularly ordered}, or $c$-ordered for short, if
there exists an integer $h$ with $1\le h\le n$ such that
\begin{equation}\label{eq1}
0\le x_h \le x_{h+1}\le\cdots\le x_n\le x_1\le\cdots\le x_{h-1}< 1
\end{equation}
If (\ref{eq1}) holds, then either $x_1=\cdots=x_n$, or the integer $h$
 is unique. Also, if $(x_1,\ldots, x_n)$ is $c$-ordered,
then clearly $(x_2,\ldots, x_n, x_1)$ is $c$-ordered. Any subsequence of a
$c$-ordered sequence is $c$-ordered. Observe also that if
$(x_1,\ldots, x_n)$ is $c$-ordered and $x_1<x_n$ then $x_1\le\cdots
\le x_n$. Indeed, if (\ref{eq1}) holds for $h\ne1$, then $x_n\le x_1$.

Two rules are useful.
\begin{lemma} 
{\rm(i)\ Translation Rule} If  $(x_1,\ldots, x_n)$ is $c$-ordered and
$y_i\equiv x_i+\alpha\bmod1$, then  
$(y_1,\ldots, y_n)$ is $c$-ordered. 

\noindent{\rm(ii)\ Insertion Rule} If  $(x_1,\ldots, x_n)$ and $(y_1,\ldots, y_m)$ are $c$-ordered, if furthermore $y_1\ne y_m$ 
and $x_i=y_1$, $x_{i+1}=y_m$, then 
$(x_1,\ldots,x_i,y_2, \ldots,y_{m-1}, x_{i+1}, \ldots ,x_n)$
is $c$-ordered.
\end{lemma}
\begin{proof}.\ (i) We may assume $0\le x_1\le \cdots \le  x_n < 1$. The
real numbers $x_i+\alpha$ satisfy
$x_1+\alpha\le \cdots \le x_n+\alpha <1+ x_1 +\alpha$. If $x_n+\alpha < 1$,
then $y_i= x_i+\alpha$ and $(y_1,\ldots, y_n)$ is $c$-ordered.
Otherwise, let $h$ be the smallest integer such that $x_h+\alpha \ge
1$. Then
$$x_1+\alpha \le  \cdots \le  x_{h-1} + \alpha < 1 \le x_h+\alpha\le \cdots \le 
x_n+\alpha$$
If $h=1$, one gets $1<x_1+\alpha \le  \cdots \le x_n+\alpha <2$ and clearly
$(y_1,\ldots, y_n)$ is $c$-ordered. If $h>1$, then $x_n+\alpha -1 <
x_1 + \alpha$ implies
$$y_h\le \cdots \le y_n<y_1\le \cdots\le y_{h-1}$$
(ii) There are two cases. If $x_i=\max\{x_1,\ldots,x_n\}$, then
$x_{i+1}\le \cdots\le x_n\le x_1\le \cdots\le x_i$. From $x_i=y_1$, $x_{i+1}=y_m$, it
follows that $y_m < y_1$. Let $h\ne1$ be the integer such that
$0\le y_h\le \cdots\le y_m<y_1\le \cdots\le y_{h-1}$. Then
$$0\le
y_h\le \cdots\le y_m=x_{i+1}\le \cdots\le x_n\le x_1\le \cdots\le x_i=y_1\le \cdots\le y_{h-1}$$
If $x_i<\max\{x_1,\ldots,x_n\}$, then $x_i=y_1<y_m=x_{i+1}$ and
consequently $x_i=y_1\le y_2\le \cdots\le y_m=x_{i+1}$.
\end{proof}
We observe that the insertion rule does not hold if $y_1=y_m$. Consider
the two $c$-ordered sequences $(x,x,y)$ and $(x,y,x)$, where $0<x<y<1$. Inserting the second into the first give the sequence $(x,y,x,y)$ which is not $c$-ordered.

We prove another useful formula.
\begin{lemma}\label{lemme1} Let $\alpha<1/2$.
If $(x,y,x+\alpha)$ is $c$-ordered, then 
$(x,y,x+\alpha,y+\alpha)$ is $c$-ordered.
\end{lemma}
\begin{proof}.\ The condition $\alpha<1/2$ implies that $(x,x+\alpha,
x+2\alpha)$ is $c$-ordered. By the translation rule, we get that 
$(x+\alpha, y+\alpha, x+2\alpha)$ is $c$-ordered. The insertion rule
shows that $(x,x+\alpha, y+\alpha, x+2\alpha)$ is $c$-ordered and,
again by the insertion rule, one gets that 
$(x,y,x+\alpha,y+\alpha)$ is $c$-ordered.
\end{proof}
\section{Intersection}
Circular order is useful in considering intersections of intervals.
Let $I = [x,y[$ be a nonempty interval. Then $x'\in [x,y[$ iff $(x,x',y)$ is
ordered.
Let $I = [x,y[$ and $I'=[x',y'[$ be nonempty intervals. Then $x'\in I$ iff $(x,x',y)$ is
$c$-ordered. Since $I\cap
I'\ne\emptyset$ iff $x'\in I$ or $x\in I'$, the intervals $I$ and $I'$
are disjoint iff $(x,y,x')$ and $(x',y',x)$ are
$c$-ordered. Consequently,
we have shown
\begin{lemma} \label{lemme3}
Let $I = [x,y[$ and $I'=[x',y'[$ be nonempty intervals. Then $I\cap
I'=\emptyset$ if and only if $(x,y,x',y')$ is $c$-ordered. 
\QED
\end{lemma}
The \textit{length} $s$ of an interval $I=[x,y[$ is the number $s=y-x$ if
$x\le y$, and is $s= 1-(x-y)$ if $y<x$. In both cases, $y\equiv
x+s\bmod1$ so that, knowing the length, we may write $I=[x,x+s[$.
\begin{lemma} \label{lemme4}
Let $I = [x,y[$ and $I'=[x',y'[$ be intervals of the same length 
$0<\alpha<1/2$. If $I$ and $I'$ intersect, then
$(x,x',y,y')$ or $(x',x,y',y)$ is $c$-ordered. In the first case,
$I\cap I'= [x',y[$, in the second case, $I\cap I'= [x,y'[$.
\end{lemma}
Observe that if the length of $I$ and $I'$ is greater than $1/2$, 
then the intersection needs not to be an interval.
\begin{proof}.\ The discussion before Lemma~\ref{lemme3} shows
that $I$ and $I'$ intersect if and only if 
$(x,x',y)$ or $(x', x, y')$ are
$c$-ordered. From Lemma~\ref{lemme1}, it follows that 
$(x,x',y,y')$ or $(x',x,y',y)$ is $c$-ordered. Moreover, $y\ne x'$ and
$x\ne y'$ since otherwise $(x,y,x',y')$ is $c$-ordered and the intervals are disjoint by Lemma~\ref{lemme1}. If $x=x'$ (or equivalently if $y=y'$), then $I=I'$.
Thus, we may assume that the numbers $x,y,x',y'$ are distinct.

Assume the first
ordering holds. The formula for the intersection is straightforward if
$0\le x<x'<y<y'<1$. If $0\le x'<y<y'<x<1$, then $I=[0,y[\cup[x,1[$
and $I\cap I'= [x',y[$. The two other cases are proved in the same way.
\end{proof}
The previous lemma will be applied to the intervals $I_k=[\beta_k,
\beta_k+\alpha[ $. They all have same length $\alpha$. We write the
conclusion for further reference.
\begin{lemma} \label{lemme2} Let $\alpha < 1/2$.
Let $I_k=[\beta_k, \beta_k+\alpha[ $ and $I_\ell=[\beta_\ell,
\beta_\ell+\alpha [$ be two intervals. 
If $I_k$ and $I_\ell$
intersect then
$(\beta_k, \beta_\ell, \beta_k+\alpha, \beta_\ell+\alpha)$
or 
$(\beta_\ell, \beta_k, \beta_\ell+\alpha,\beta_k+\alpha)$
is $c$-ordered. Moreover, $I_k\cap I_\ell= [\beta_\ell, \beta_k+\alpha[$ in the
first case, and $I_k\cap I_\ell= [\beta_k, \beta_\ell+\alpha[$ in the
second case. \QED
\end{lemma} 

The following observation is the basic step for analyzing the coding
induced by a rotation. Recall that for $K\subset\{0,\ldots, m\}$,
$$X_K = \bigcap_{k\in K} I_k\cap \bigcap_{k\notin K} \overline{I_k}$$
We assume from now on that $\alpha<1/2$.
\begin{proposition}
Assume $X_K\ne\emptyset$ for some $K\subset\{0,\ldots, m\}$ and assume
$(\beta_{i_1}, \beta_{i_2}, \beta_{i_3}, \beta_{i_4})$ is a
$c$-ordered sequence. If $i_1, i_3\in K$, then $i_2\in K$ or $i_4\in K$. 
\end{proposition}
\begin{proof}.\ Arguing by contradiction, suppose that $i_2,i_4\notin K$. 
Since $X_K\ne\emptyset$, the interval
$I_{i_1}\cap I_{i_3}$ is not empty, therefore by Lemma~\ref{lemme2}
$(\beta_{i_1}, \beta_{i_3}, \beta_{i_1}+\alpha, \beta_{i_3}+\alpha)$
or
$(\beta_{i_3}, \beta_{i_1}, \beta_{i_3}+\alpha, \beta_{i_1}+\alpha)$
is $c$-ordered (or the sequence obtained by exchanging $i_1$ and
$i_3$). Consider the first case, the second is the same by
exchanging $i_2$ and $i_4$. 
Since 
$(\beta_{i_1}, \beta_{i_2}, \beta_{i_3})$ 
is $c$-ordered,  the translation rule shows that 
$(\beta_{i_1}+\alpha , \beta_{i_2}+\alpha , \beta_{i_3}+\alpha )$ 
is $c$-ordered which gives, applying twice the insertion rule, that
$(\beta_{i_1}, \beta_{i_2},\beta_{i_3}, 
\beta_{i_1}+\alpha, \beta_{i_2}+\alpha, \beta_{i_3}+\alpha)$
is $c$-ordered.
{}From this, we get that 
$(\beta_{i_1}, \beta_{i_2},
\beta_{i_1}+\alpha, \beta_{i_2}+\alpha)$
is $c$-ordered. 
{}From Lemma~\ref{lemme3}, we know that
$I_{i_1}\cap I_{i_3} = [\beta_{i_3},\beta_{i_1}+\alpha[$, and this
is then disjoint from 
$\overline{I_{i_2}}= [\beta_{i_2}+\alpha,\beta_{i_2}[$.
\end{proof}

\begin{proposition} If $X_K$ is not empty, then there exist integers
$k,\ell$ with $0\le k< \ell\le m$ such that $\{K,M\setminus K\}= \{
\{k,\ldots, \ell-1\}, \{\ell,\ldots m, 0,\ldots, k-1\} \}$
\end{proposition}
\begin{proof}. This is a direct consequence of the preceding discussion.
\end{proof}
It follows that there are only $(m+1)(m+2)$ intervals $X_K$ to be
considered. In fact, consider the numbers $0,\beta_1,\ldots,\beta_m,1$
and $\alpha, \beta_1+\alpha,\ldots,\beta_m+\alpha$. They partition
$[0,1[$ into exactly $2m+2$ intervals. Each of these intervals is
contained in one and only one of the $X_K$ (but $X_\emptyset$ may be
scattered over several of the small intervals). This means that, among
the $(m+1)(m+2)$ possible intervals $X_K$, there are only $2m+2$ that
are used in a particular setting of the values of $\alpha,
\beta_1,\ldots,\beta_m$.

\begin{theorem} Assume $K\ne\emptyset, M$, and $X_K\ne\emptyset$.
Then $X_K=[\beta_{\ell-1},\beta_k+\alpha[
\cap[\beta_{k-1}+\alpha,\beta_\ell[$.
\end{theorem}
If $K=\{k\}$ is a singleton, then the formula still holds with $\ell-1=k$.
\begin{proof}. 
Suppose that $K=\{k, \cdots, \ell-1 \}$ with $k < \ell$. The other
case is symmetric. We first prove that $\bigcap_{n \in K} I_n=  [\beta_{\ell
-1},\beta_k +\alpha [.$ Set $Y_K=\bigcap_{n \in K} I_n$.

Since $X_K \neq \emptyset$, the interval $I_k \cap I_{\ell-1}$ is not empty.
By Lemma 4.3 there are two cases:
 either $(\beta_k, \beta_{\ell -1},\beta_k +\alpha, \beta_{\ell -1}+\alpha)$
is $c$-ordered,
 or $(\beta_{\ell -1},\beta_k, \beta_{\ell -1}+\alpha,\beta_k +\alpha)$
is $c$-ordered.

 We show that this second case cannot happen. Indeed in this case, $Y_K \subset I_k \cap I_{\ell-1}= [ \beta_k,
\beta_{\ell-1} + \alpha[.$
Moreover for each $n \in M \setminus K$, the sequence $( \beta_{\ell-1},
\beta_n,\beta_k)$ is $c$-ordered.
By translation and insertion, the sequence   $( \beta_{\ell-1}, \beta_n,\beta_k, \beta_{\ell-1} +
\alpha, \beta_n + \alpha,\beta_k + \alpha)$ is $c$-ordered. This shows
that $I_n \supset I_k \cap I_{\ell-1}$, and consequently 
$I_k \cap I_{\ell-1} \cap \overline{I_n} = \emptyset$  for each  $n$ in $M \setminus K$,
 contradicting the assumption that $X_K \neq \emptyset.$

Thus, $(\beta_k, \beta_{\ell -1},\beta_k +\alpha, \beta_{\ell -1}+\alpha)$
is $c$-ordered. This   implies that $I_k \cap I_{\ell-1}= [\beta_{\ell -1},\beta_k
+\alpha [.$
If $i \in K$ then $(\beta_k, \beta_i,\beta_{\ell -1})$ is $c$-ordered.
By translation, $(\beta_k + \alpha , \beta_i + \alpha ,\beta_{\ell -1} + \alpha)$ is $c$-ordered.
By insertion of  $(\beta_k, \beta_i,\beta_{\ell -1})$  into $(\beta_k,
\beta_{\ell -1},\beta_k +\alpha, \beta_{\ell -1}+\alpha)$
one gets  $(\beta_k, \beta_i, \beta_{\ell -1},\beta_k +\alpha
,\beta_{\ell -1}+\alpha)$ is $c$-ordered.
Again by insertion of $(\beta_k + \alpha , \beta_i + \alpha
,\beta_{\ell -1} + \alpha)$,
the sequence 
 $(\beta_k, \beta_i, \beta_{\ell -1},\beta_k +\alpha,
\beta_i+\alpha ,\beta_{\ell -1}+\alpha)$ is $c$-ordered. Thus $Y_K=  [\beta_{\ell -1},\beta_k +\alpha [.$

The second part of the proof deals with $\bigcap_{n \in M \setminus K}
\overline{I_n}.$ In this intersection the index $n$ runs through the set $\{0, \cdots, k-1,\ell, \cdots, m\}$.
The set $M \setminus K$ is partitioned into three possibly empty
subsets as follows: $n \in N$ iff $ \overline{I_n} \supset Y_K$,
 $n \in P$  iff
 $ \overline{I_n} \cap Y_K = [\beta_{\ell-1}, \beta_n[$
and finally $n \in Q$ iff
$ \overline{I_n} \cap Y_K= [ \beta_n +\alpha, \beta_k
+ \alpha[.$ Of course, $$
X_K= \bigcap_{n \in N}  (\overline{I_n} \cap Y_K) \cap
\bigcap_{n \in P}  (\overline{I_n}  \cap Y_K) \cap  \bigcap_{n \in Q}
(\overline{I_n} \cap Y_K)$$
If one of the sets $N,P,Q$ is empty it does not contribute to the intersection.

Clearly $ \bigcap_{n \in N}  (\overline{I_n} \cap Y_K)= Y_K$.
Next $\bigcap_{n \in P}  (\overline{I_n}  \cap Y_K)= \bigcap_{n \in P}
[\beta_{\ell-1}, \beta_n[$.
If $P$ is not empty then $\ell$ is in $P$ and  $\bigcap_{n \in P}  (\overline{I_n}  \cap Y_K)= 
[\beta_{\ell-1}, \beta_\ell[$. 
Finally,  $\bigcap_{n \in Q}  (\overline{I_n}  \cap Y_K)= \bigcap_{n \in Q}
[\beta_{n} + \alpha, \beta_k + \alpha[$.
If $Q$ is not empty then $k-1$ is in $Q$ and  $\bigcap_{n \in Q}  (\overline{I_n}  \cap Y_K)= 
[\beta_{k-1} + \alpha , \beta_k + \alpha[$.

To finish the proof, we just have to verify that
in each case,  $X_K=\bigcap_{n \in N}  (\overline{I_n} \cap Y_K) \cap
\bigcap_{n \in P}  (\overline{I_n}  \cap Y_K) \cap  \bigcap_{n \in Q}
(\overline{I_n} \cap Y_K)$
 is equal to $[\beta_{\ell-1},\beta_k+\alpha[
\cap[\beta_{k-1}+\alpha,\beta_\ell[$.

If $P \neq \emptyset$ then  $\bigcap_{n \in P}  (\overline{I_n}  \cap Y_K)= 
[\beta_{\ell-1}, \beta_\ell[$ and the sequence $(\beta_{\ell-1}, \beta_\ell,
\beta_k + \alpha)$ is $c$-ordered (case $P_1$).
 If
 $P = \emptyset$ then  the sequence $(\beta_{\ell-1},
\beta_k + \alpha, \beta_\ell)$ is $c$-ordered (case $P_0$).
If $Q \neq \emptyset$ then  $\bigcap_{n \in Q}  (\overline{I_n}  \cap Y_K)= 
[\beta_{k-1} + \alpha , \beta_k + \alpha[$  and the sequence
$(\beta_{\ell-1}, \beta_{k-1} + \alpha,
\beta_k + \alpha)$ is $c$-ordered (case $Q_1$).
 If
 $Q = \emptyset$ then  the sequence $(\beta_{\ell-1},
\beta_k + \alpha, \beta_{k-1} + \alpha)$  is $c$-ordered (case $Q_0$).

Case $(P_1 Q_1)$. If $P$ and $Q$ are nonempty then 
$X_K=[\beta_{\ell -1},\beta_k +\alpha [ \cap
[\beta_{\ell-1}, \beta_\ell[ \cap [\beta_{k-1} + \alpha, \beta_k +
\alpha[$. As the sequences $(\beta_{\ell-1}, \beta_\ell,
\beta_k + \alpha)$ and $(\beta_{\ell-1}, \beta_{k-1} + \alpha,
\beta_k + \alpha)$ are $c$-ordered, by the  insertion rule
either  the sequence
$(\beta_{\ell-1}, \beta_\ell, \beta_{k-1} + \alpha, \beta_k + \alpha)$ or  
$(\beta_{\ell-1},\beta_{k-1} + \alpha,   \beta_\ell,  \beta_k +
\alpha)$ is $c$-ordered.
The first case is impossible because $X_K$ is not empty.
The second case implies that $X_K=[\beta_{\ell-1},\beta_k+\alpha[
\cap[\beta_{k-1}+\alpha,\beta_\ell[$.

Case $(P_0 Q_1)$. If $P = \emptyset$ and $Q \neq  \emptyset$  then 
$X_K= [\beta_{\ell-1},\beta_k+\alpha[
\cap[\beta_{k-1}+\alpha,\beta_k + \alpha[$ and
the sequences $(\beta_{\ell-1},
\beta_k + \alpha, \beta_\ell)$,  $(\beta_{\ell-1}, \beta_{k-1} + \alpha,
\beta_k + \alpha)$ are $c$-ordered. By insertion 
 the sequence $(\beta_{\ell-1},\beta_{k-1} + \alpha,
\beta_k + \alpha, \beta_\ell)$ is $c$-ordered. Thus
$X_K=[\beta_{\ell-1},\beta_k+\alpha[
\cap[\beta_{k-1}+\alpha,\beta_\ell[$.

Case $(P_1 Q_0)$ is symmetric to case  $(P_0 Q_1)$.

Case $(P_0 Q_0)$. If $P = \emptyset$ and $Q =  \emptyset$  then 
$X_K= [\beta_{\ell-1},\beta_k+\alpha[
\cap[\beta_{k-1}+\alpha,\beta_\ell[$ and
the sequences $(
\beta_k + \alpha,\beta_{k-1} + \alpha, \beta_{\ell -1})$,
  $( \beta_{k-1} + \alpha,
\beta_\ell,\beta_{\ell-1})$ are $c$-ordered. By insertion rule
either  the sequence
$(\beta_k + \alpha,\beta_{k-1} + \alpha, \beta_\ell, \beta_{\ell -1})$, or  
$(\beta_k + \alpha,\beta_{k-1} + \alpha, \beta_\ell, \beta_{\ell -1})$
is $c$-ordered.
The first case is impossible because $X_K$ is non empty.
The second case implies that $X_K=[\beta_{\ell-1},\beta_k+\alpha[
\cap[\beta_{k-1}+\alpha,\beta_\ell[$.
\end{proof}

Remark: As an additional property, the preceding proof shows that $X_K$ is an
interval and the interior of $X_K$
does not contain any $\beta_i$ or $\{ \beta_i + \alpha\}.$

\section{Main result}

\begin{proposition}
If $x\in B_i$ and $x+\alpha\in B_j\cap X_K$, then $j\equiv i+ |K|\bmod m+1$.
\end{proposition}

\begin{proof}
If $K \neq \emptyset$ or $K \neq M$ then 
$X_K= [\beta_{\ell -1},\beta_k +\alpha [ \cap [\beta_{k-1}+\alpha,
\beta_\ell[.$
As $x \in [ \beta_i, \beta_{i+1}[,$ by translation rule we have $ y = x+
\alpha \in [ \beta_i + \alpha, 
\beta_{i+1}+\alpha[$. Furthermore, $y \in [\beta_j, \beta_{j+1}[.$
By the preceding remark and by identification, the only possibility is $j=\ell -1$ and $i=
k-1.$
It follows that $j=i + |K| \bmod m+1.$

If $K= \emptyset$ then $X_K= \cap_{n \in M}
\overline{I_n}.$
By hypothesis we have $y=x+ \alpha \in B_j.$ But $y \in X_K$ implies
that $y \notin [\beta_j, \beta_j + \alpha[ = I_j 
.$

If $|I_j| \geq |B_j|$ then $B_j \cap X_K$ should be empty in
contradiction with the hypothesis.
Thus $|I_j| < |B_j|$. The interval $B_j$ is equal to $I_j \cup [\beta_j
+ \alpha, \beta_{j+1}[.$
That is $x+ \alpha \in  [\beta_j + \alpha, \beta_{j+1}[$ and $x+ \alpha
\in [\beta_i + \alpha, \beta_{i+1} + \alpha[.$
By identification, we have $i=j.$

If $K=M$ then $X_K= \cap_{n \in M} I_n.$
As $x \in I_n$ implies $x+\alpha \notin I_n$, by contraposition $x+
\alpha \in  X_K$ 
implies $x \notin I_n$ for all  $n$. Thus  $(x,\beta_n,x+\alpha)$ is
$c$-ordered for all $n \in M$. As $x \in  B_i$ and $x$ is not in $I_i$
the sequence $(\beta_i, \beta_i + \alpha, x, \beta_{i+1})$ is $c$-ordered.
Thus  $(x, \beta_{i+1}, \cdots , \beta_{i+m-1}, \beta_{i}, x + \alpha)$ is
$c$-ordered.
Consequently  $X_K$ is equal to $[\beta_{i},
\beta_{i+1}+\alpha[.$
As $x+ \alpha \in [\beta_j, \beta_{j+1}[$ by identification we find
$j=i$.
\end{proof}

{}From this proposition, we get the following automaton $\A$
 (Figure \ref{Tr} gives the automata for  $m=1,2,3$). 
Its set of states is the set $M$ in bijection with the intervals $B_k$. 
The alphabet is the set of subsets of $M$ corresponding to the
nonempty intervals $N_K$. As already mentioned, there are $(m+1)(m+2)$ of them.
The transitions or edges are given by the proposition: $(i, K, j)$ is
an transition if $j\equiv i+ |K|\bmod m+1$.

Observe that the automaton is deterministic. Also, it is universal in
the following sense : for a particular setting of $\alpha,
\beta_1,\ldots,\beta_m$, if the $\beta_i$ and the $\beta_j + \alpha$
are two by two distinct, there are only $2m+2$ of the edges that are
used. Indeed they are exactly $2m+2$ intervals in the partition
and  between $\beta_i$ and $\beta_{i+1}$ the coding is uniquely
determined for all $i$.

\end{document}